\documentclass[a4paper]{article}
\usepackage{amsmath, amstext, amsxtra, amsthm, amscd, amsgen,amsbsy, amsopn, amsfonts, latexsym, amssymb, amscd, amsthm}
\usepackage{marvosym}
\usepackage{color}
\usepackage{todonotes}
\usepackage{tikz-cd}

\def\proofend{\hbox to 1em{\hss}\hfill $\blacksquare $\bigskip }

\newtheorem*{theorem*}{Main Theorem}

\newtheorem{theorem}{Theorem}[section]
\newtheorem{proposition}[theorem]{Proposition}
\newtheorem{lemma}[theorem]{Lemma}

\def\Z{{\mathbb Z}}
\def\R{{\mathbb R}}

\def\C{{\mathbb C}}
\def\N{{\mathbb N}}

\begin{document}


\title{Moduli space of metrics of nonnegative sectional or positive Ricci curvature on Brieskorn quotients in dimension 5}

\author{Jonathan Wermelinger}

\date{\today}

\maketitle
\begin{abstract}\noindent We show that the moduli space of Riemannian metrics of nonnegative sectional curvature on certain quotients of Brieskorn varieties in dimension 5 has infinitely many path components. The same is true for the corresponding moduli space of positive Ricci curvature.\end{abstract}

\section{Introduction}
The subject of moduli spaces of Riemannian metrics with lower curvature bounds has been an active field of study over the last years. The first results on the moduli spaces of 5-manifolds have been presented in \cite{TWi17}. Recently, the first examples of 5-manifolds whose moduli spaces have infinitely many path-components for nonnegative sectional curvature \cite{DGA19} and positive Ricci curvature \cite{G20} have been identified.

In this note, we continue the study of moduli spaces of metrics of nonnegative sectional curvature of quotients of Brieskorn varieties in dimension 5 which was started by Dessai and Gonz\'{a}lez-\'{A}lvaro. One of their results is about the moduli spaces of manifolds homotopy equivalent to $\mathbb{R}\text{P}^5$ whose universal cover is $S^5$, while we will consider quotients of $S^2\times S^3$ under certain involutions. This really only changes the diffeomorphism classification of or quotients and we may apply a very similar procedure as presented in \cite{DGA19}. We also want to point out that our manifolds are also different from the examples presented by Goodman \cite{G20}, since his manifold's second homotopy group is trivial as a $\mathbb{Z}[Z_2]$-module and our's are non-trivial as we will discuss below.





The main result we are going to prove is the following.

\begin{theorem*}
Let $N^5$ be an orientable, closed, smooth $5$-dimensional manifold whose universal cover is $S^2\times S^3$, satisfying $\pi_1(N)=\Z_2$, $w_2(N)\neq 0$ and $\pi_2(N)$ is non-trivial as a $\Z[\Z_2]$-module. Then the moduli space $\mathcal{M}_{sec\geq 0}(N)$ of nonnegative sectional curvature has infinitely many path components. The same is true for the moduli space $\mathcal{M}_{Ric> 0}(N)$ of metrics of positive Ricci curvature.
\end{theorem*}

As we will see, there are five pairwise non-diffeomorphic manifolds satisfying the conditions of this theorem. These are all represented by what we will call a Brieskorn quotient. Furthermore, these manifolds are not all pairwise homotopy equivalent. Unless otherwise stated, all manifolds will be oriented if orientable and all diffeomorphisms will be orientation-preserving.

This note is organized as follows. In section 2, we define Brieskorn manifolds and a cohomogeneity one action on them. From this action, we select an involution and define Brieskorn quotients. We discuss the diffeomorphism classification of manifolds satisfying the assumptions of the main theorem and relate them to the Brieskorn quotients. In section 3, we equip the quotients and other constructions with Riemannian metrics of nonnegative sectional curvature. Section 4 then introduces $Spin^c$-structures and eta-invariants. Finally, the corresponding moduli space is investigated and the proof of our main theorem is given in section 5.

\section{Brieskorn varieties and manifolds}

\subsection{Definition}\label{subsec:def}

Let $z=(z_0,z_1,z_2,z_3)\in\C^4$ and consider the polynomial $$f:\C^{4}\rightarrow \C: z\mapsto z_0^d+z_1^2+z_2^2+z_3^2.$$
For $\epsilon \in \C$ we get the hypersurface $V_\epsilon^3(d):=f^{-1}(\epsilon)$, which is an algebraic manifold with no singular points for $\epsilon\neq 0$, but only a variety if $\epsilon=0$ with an isolated singular point $z=0$ if $d\geq 2$. Let $S^{7}:=\{z\in \C^{4}|\text{ }|z_0|^2+...+|z_3|^2=1\}$ and $D^{8}:=\{z\in \C^{4}|\text{ }|z_0|^2+...+|z_3|^2\leq 1\}$ denote the unit sphere and the closed disk in $\C^{4}=\R^{8}$ respectively. We may then define the link $$\Sigma^{5}_\epsilon(d):=V_\epsilon^3(d)\cap S^{7}$$ and $$W^{6}_\epsilon(d):=V_\epsilon^3(d)\cap D^{8},$$ both equipped with the natural orientation. The manifolds $\Sigma^5_\epsilon(d)$ are examples of \emph{Brieskorn varieties} (see \cite{HM68} and \cite{M68} for more details).

An important fact is that for $\epsilon$ small enough, $\Sigma^{5}_\epsilon(d)$ is diffeomorphic to $\Sigma^{5}_0(d)$ (see \cite[Satz 14.3]{HM68}). Also, $W^{6}_\epsilon(d)$ is homotopy equivalent to a bouquet $S^3\vee ... \vee S^3$ of $3$-spheres (see \cite[p.83]{HM68}).

If $d$ is odd, then $\Sigma^5_0(d)$ is diffeomorphic to $S^5$ (see \cite[14.8 (zweiter) Satz]{HM68}).

Throughout this note, we will consider the case when $d$ is even. In this case, it can be shown that $\Sigma^5_0(d)$ is diffeomorphic to $S^2\times S^3$ (see \cite[Proposition 7]{GT98}).

\subsection{Group action and quotient}\label{subsec:grpac}

Recall that an action of a compact Lie group $G$ on a smooth manifold $M$ is said to be of \emph{cohomogeneity one} if the orbit space $M/G$ is one-dimensional. Consider the following group actions. $$a:S^1\times \C^4 \rightarrow \C^4: (w,(z_0,z_1,z_2,z_3))\mapsto (w^2z_0,w^dz_1,w^dz_2,w^dz_3),$$ $$b:O(3)\times \C^4 \rightarrow \C^4: (A,(z_0,z_1,z_2,z_3))\mapsto (z_0,(A(z_1,z_2,z_3)^T)^T).$$

Suppose $d$ is even. In case $\epsilon\neq 0$, we get an action of $\Z_{2d}\times O(3)$ with ineffective kernel $\{(\pm 1, Id)\}$, whereas for $\epsilon = 0$ we get an action of $S^1\times O(3)$ with the same ineffective kernel as before. Both of these actions are orientation-preserving on $\Sigma^{5}_\epsilon(d)$.

The action of $S^1\times O(3)$ on $\Sigma^{5}_0(d)$ is of cohomogeneity one with orbit space an interval. This can be seen by noticing that the norm of the zero component is not affected by this action.

The principal isotropy is $\Z_2\times O(1)$, corresponding to interior points of the quotient interval. The singular isotropies are $S^1\times O(1)$ and $\Z_2\times O(2)$, corresponding to the vanishing and maximal norm of $z_0$ respectively. The corresponding singular orbits are both of codimension two. For more details on the above cohomogeneity one action, see \cite[\textsection 4.2]{DGA19} and \cite[\textsection 1]{GVWZ06}.

Consider $\tau:=(1,-Id)\in S^1\times O(3)$, i.e. $\tau(z)=(z_0,-z_1,-z_2,-z_3)$ for $z=(z_0,z_1,z_2,z_3)\in\C^4$. This gives an orientation-preserving involution which acts freely on $\Sigma^{5}_\epsilon(d)$ for $\epsilon<1$. Since the polynomial $f(z)=z^d_0+z^2_1+z^2_2+z^2_3$ is $\tau$-invariant, we get an induced map $\tilde{f}:S^{7}/\tau\rightarrow \C$.

We can now define $N^5_\epsilon(d):=\Sigma^{5}_\epsilon(d)/\tau$, which we will call a \emph{Brieskorn quotient}. For $d$ odd, these are homotopy equivalent to $\R\text{P}^5$. For $d$ even, these quotients are not all pairwise homotopy equivalent (see \cite[p.1180 Remark (2)]{GT98}). In both cases, $N^5_\epsilon(d)$ is diffeomorphic to $N^5_0(d)$ for $\epsilon$ sufficiently small (see \cite[Proposition 4.1]{DGA19}).

The above action of $S^1\times O(3)$ on $\Sigma^{5}_0(d)$ descends to a cohomogeneity one action on $N^5_\epsilon(d)$ which also has codimension two singular orbits.

For $0<\epsilon<1$, the fixed points of the action of $\tau$ on $W^{6}_\epsilon(d)$ are $p_i=(\lambda_i,0,0,0)$ , $1\leq i\leq d$, where $\lambda_i$ is a complex d-root of $\epsilon$ for all $i$. These fixed points are isolated and all lie in the interior of $W^{6}_\epsilon(d)$. The subaction $\Z_{2d}\times \{Id\}$ of $\Z_{2d}\times O(2)$ on $W^{6}_\epsilon(d)$ permutes these fixed points.

\subsection{Diffeomorphism classification}\label{subsec:diffclass}

In order to give a diffeomorphism classification of the above quotients, we need to introduce $Pin$ structures. Define the groups $Pin^\pm(n)$ to be the two distinct double covers of $O(n)$, corresponding to the central extensions of $O(n)$ by $\Z_2$: $$1\rightarrow \Z_2\rightarrow Pin^\pm(n)\rightarrow O(n) \rightarrow 1.$$
They are topologically the same space but have different group structures. A $Pin^\pm$-structure on a smooth manifold $M$ is a reduction of the structure group $O(n)$ of its tangent bundle to $Pin^\pm(n)$.

A manifold $M$ admits a $Pin^+$-structure if and only if $w_2(M)=0$ and a $Pin^-$-structure if and only if $w_2(M)=w_1(M)^2$, for $w_1,w_2$ the first and second Stiefel-Whitney class respectively. Furthermore, $Pin^\pm$-structures are in bijection with $H^1(M;\Z_2)$, if any exist (see \cite{KT90}).

If we ask for $Pin^\pm$-structures on manifolds with boundaries to restrict to $Pin^\pm$-structures on the boundary, we get $Pin^\pm$-bordism groups $\Omega^{Pin^\pm}_n$. Kirby and Taylor determined $\Omega^{Pin^+}_4\cong \Z_{16}$ (see \cite{KT90}), and if we identify an element with its inverse, we get $\Omega^{Pin^+}_4/\pm\cong \Z_{8}$.

Recall that a \emph{characteristic submanifold} of a manifold $M^n$ with $\pi_1(M)=\Z_2$ is a submanifold $P^{n-1}\subset M^n$ such that there exists a submanifold $A\subset \tilde{M}$ with $\tilde{M}=A\cup T(A)$, $\partial A = \partial T(A)=\tilde{P}$ and $P=\tilde{P}/T$. Here $\tilde{P}$ is the universal covering and $T$ the corresponding deck-transformation.

Su now proved the following classification result (see \cite[Theorem 1.2.]{S12}).

\begin{theorem}\label{thm:classification5mfds}
Let $M^5$ be an orientable, smooth 5-manifold with $\pi_1(M)\cong \Z_2$ and universal cover $\tilde{M}\cong S^2\times S^3$. If $\pi_2(M)\cong \Z$ is non-trivial as a $\Z[\Z_2]$-module and $w_2(M)\neq 0$, then $M$ is diffeomorphic to $N^5_0(d)$ for some $d=0,2,4,6,8$. The quotients $N^5_0(d)$ are classified by the $Pin^+$-bordism class of their characteristic submanifold $P_d$ with $[P_d]=d\in\Omega^{Pin^+}_4/\pm=\{0,1,...,8\}$.
\end{theorem}

The module structure on $\pi_2(M)$ comes from the $\pi_1(M)$-action on it.


Geiges and Thomas showed that for $d$ even, $N^5_0(d)$ satisfies all of the above conditions (see \cite[Proposition 6, 7 and Lemma 11]{GT98}). In particular, they show that the induced isomorphism $\tau_*:H_2(\Sigma^5_0(d);\mathbb{Z})\rightarrow H_2(\Sigma^5_0(d);\mathbb{Z})$ is non-trivial (see \cite[Lemma 5]{GT98}) and therefore it follows by the use of Hurewicz' map that $\pi_2(N^5_0(d))$ is non-trivial as a $\mathbb{Z}[\mathbb{Z}_2]$-module. Note furthermore that the $Pin^\pm$-bordism class does not depend on the choice of a characteristic submanifold (see \cite[Lemma 9]{GT98}).

From the Proof of \cite[Lemma 11]{GT98} and using the above theorem, one can also deduce the following.

\begin{lemma}\label{lemma:diffquo}
For $d,d'\in\N$ even, $N^5_\epsilon(d)$ and $N^5_\epsilon(d')$ are diffeomorphic if and only if $d\equiv \pm d' \text{ mod }16$.
\end{lemma}

We thus get five different oriented diffeomorphism types for the quotients. For a fixed $d$ even, and hence a fixed diffeomorphism type, there are infinitely many $d'\in \N$ such that $N^5_\epsilon(d)$ is diffeomorphic to $N^5_\epsilon(d')$.

\section{Metrics of nonnegative sectional curvature}\label{sec:metrics}

By the work of Grove and Ziller, one can equip cohomogeneity one manifolds with an invariant metric of nonnegative sectional curvature, provided the singular orbits are of codimension two (see \cite[Theorem E]{GZ00}).

As we have seen above, the quotients $N^5_0(d)$ admit a cohomogeneity one action by $S^1\times O(3)$ with codimension two singular orbits. Hence, $N^5_0(d)$ can be equipped with an invariant metric of $sec\geq 0$, which we will denote by $g^d_{GZ}$ and call its Grove-Ziller metric. 

In order to get a metric of positive scalar curvature (and $sec\geq 0$), we will apply Cheeger deformation to the Grove-Ziller metric. See \cite[\textsection 6.1]{AB15} for a general treatment of Cheeger deformation.

Let $G$ be a compact Lie group acting by isometries on a Riemannian manifold $(M,g)$ and $Q$ a biinvariant metric on $G$. The action $g_1\cdot(p,g_2)=(g_1\cdot p,g_1\cdot g_2)$ of $G$ on $M\times G$ is free and isometric with respect to $g+\frac{1}{t}Q$, where $t\in\R^{>0}$. The quotient $(M\times G)/G$ is diffeomorphic to $M$. The \emph{Cheeger deformation} of $g$ is the metric $g_t$ which turns the projection $$\pi_t:(M\times G,g+\frac{1}{t}Q)\rightarrow (M,g_t):(p,u)\mapsto u^{-1}p$$ into a Riemannian submersion. The Lie group $G$ acts by isometries on $(M,g_t)$ for all $t>0$. It can be shown that $g_t$ varies smoothly in $t$ and that it extends smoothly to $t=0$ by setting $g_0=g$.

Let $p\in M$, $G_p$ the isotropy group and $\mathfrak{g}_p$ the corresponding Lie algebra. Denote by $\mathfrak{m}_p$ the $Q$-orthogonal complement of $\mathfrak{g}_p$ in $\mathfrak{g}$, the Lie algebra of $G$. For $X\in \mathfrak{g}$ and $p\in M$, let $X^*_p=\frac{d}{dt}(exp(tX)\cdot p)|_{t=0}$. One can then identify the orthogonal complement with the tangent space of the $G$-orbit of $p$, $$\mathfrak{m}_p\xrightarrow{\cong}T_p(Gp):X\mapsto X^*.$$ Consider the automorphism $P_p:\mathfrak{m}_p\rightarrow \mathfrak{m}_p$ defined by the relation $Q(P_p(X),Y)=g(X^*_p,Y^*_p)$ for all $Y\in \mathfrak{m}_p$.

Cheeger deformations preserve nonnegative sectional curvature and under certain circumstances, also produce positive scalar curvature. For $X,Y\in\mathfrak{g}$, let $[X,Y]$ denote the Lie bracket of $X$ and $Y$. Then one gets the following (see \cite[Proposition 5.2]{DGA19}).

\begin{proposition}\label{Prop:cheegposscal}
Let $(M,g)$ be a Riemannian manifold and $G$ a compact Lie group acting on it by isometries. If $sec_g\geq 0$, then $sec_{g_t}\geq 0$ for all $t\geq 0$. If in addition there exist $X,Y\in \mathfrak{m}_p$ such that $[P_p(X),P_p(Y)]\neq 0$ for all $p\in M$, then $scal_{g_t}>0$ for all $t>0$.
\end{proposition}

Let us apply this to $(N^5_0(d),g^d_{GZ})$. The action is by $S^1\times O(3)=SO(2)\times O(3)$ with principal isotropy group $\Z_2\times O(1)$ and singular isotropy groups $S^1\times O(1)$ and $\Z_2\times O(2)$. Since $S^1\times O(3)$ is non-abelian and the Lie algebra of $\Z_2\times O(1)$ is $0$-dimensional, there exist non-commuting elements in this case. Now consider the singular isotropy groups, which are both 1-dimensional. It follows that the dimension of $\mathfrak{m}_p$ is 3 in both cases. By the classification of maximal tori of orthogonal groups, we know that $G=SO(2)\times O(3)$ has a 2-dimensional abelian maximal torus. Its Lie algebra $\mathfrak{g}$ is 4-dimensional and therefore splits into a 2-dimensional abelian and a 2-dimensional non-abelian part. $P_p$ being an automorphism, this means that there must be non-commuting elements in its image. Thus the conditions of Proposition \ref{Prop:cheegposscal} are satisfied.

Let $t'>0$ be fixed and $g^d:=(g^d_{GZ})_{t'}$ the Cheeger deformation of the Grove-Ziller metric at time $t'$. The metric $g^d$ will be called the \emph{modified Grove-Ziller metric} on $N^5_0(d)$, which is of $sec\geq 0$, $scal>0$ (by Proposition \ref{Prop:cheegposscal}) and invariant under the action of $S^1\times O(3)$.

We can now extend this metric to

$$\tilde{Z}(d):=\tilde{f}^{-1}([0,\epsilon_0])=\bigcup_{0\leq \epsilon \leq \epsilon_0}N^5_\epsilon(d),$$
where $\tilde{f}:S^7/\tau\rightarrow \C$ is the map induced from the polynomial $f$ discussed in \textsection \ref{subsec:grpac}. This extended metric can further be averaged over $O(3)$, yielding an $O(3)$-invariant metric $h$ on $\tilde{Z}(d)$.

If we restrict $h$ to $N^5_{\epsilon_0}(d)$, we can lift it to $\Sigma^5_{\epsilon_0}(d)$ via pull-back and extend it to a metric on $W^6_{\epsilon_0}(d)$ in such a way that it is of product form near the boundary. Furthermore, the latter metric can be assumed to be $O(3)$-invariant and of positive sectional curvature when restricted to a sufficiently small invariant neighborhood of each fixed point of the action of $\tau$ on $W^6_{\epsilon_0}(d)$ (see \cite{DGA19} Corollary 5.5 and discussion thereafter). After averaging, we finally get an $O(3)$-invariant metric $k$ on $W^6_{\epsilon_0}(d)$.

Now, using $O(3)$-actions, we can produce Cheeger deformations $g^d_t$, $h_t$ and $k_t$ of the above metrics.  These metrics satisfy the following properties (see \cite[\textsection 5 and 6]{DGA19}).

\begin{lemma}\label{lemma:metrics}
Let $g^d_t$, $g_t$ and $k_t$ be the Cheeger deformations of the metrics $g^d$, $h$ and $k$ from above. Then,

\begin{itemize}
\item for every $t> 0$, the metric $g^d_t$ is of $sec\geq 0$ and $scal>0$,
\item there exists a $t_0\geq 0$ such that the restriction of $h_t$ to $N^5_{\epsilon}(d)$ has $scal>0$ for all $t\geq t_0$ and $0\leq \epsilon\leq \epsilon_0$, and
\item $k_t$ on $W^6_{\epsilon_0}(d)$ also has $scal>0$ and is of product form near the boundary for all $t\geq t_0$. 
\end{itemize}
\end{lemma}

\section{$Spin^c$-structures and reduced eta-invariants}

The eta-invariant we want to use to distinguish path components in the space of metrics requires some additional structure on the considered manifolds. In this section, we give a short introduction on these $Spin^c$-structures (see \cite[Appendix D]{LM89} for more details).

Let $Spin^c(n):=Spin(n)\times_{\Z_2}U(1)$, where the action is diagonal by $\{\pm(1,1)\}$. Then there is a short exact sequence
$$1 \rightarrow \Z_2\rightarrow Spin^c(n)\xrightarrow{\xi} SO(n)\times U(1) \rightarrow 1$$ where $\xi$ is the map induced by $\rho\times (-)^2$ with $\rho:Spin(n)\rightarrow SO(n)$ the standard double-covering. A $Spin^c$-structure on an n-dimensional, oriented Riemannian manifold $(X,g)$ consists of a principal $Spin^c(n)$-bundle $P_{Spin^c(n)}$ over $X$, a principal $U(1)$-bundle $P_{U(1)}$ over $X$ and a bundle map $\xi:P_{Spin^c(n)}\rightarrow P_{SO(n)}\times P_{U(1)}$ such that $\xi(pg)=\xi(p)\rho(g)$ for all $p\in P_{Spin^c(n)}$ and $g\in Spin^c(n)$. Here $P_{SO(n)}$ denotes the principal bundle of oriented orthonormal frames on $X$. If $X$ is equipped with a $Spin^c$-structure, it is called a $Spin^c$-manifold.

A Riemannian manifold $X$ can be equipped with a $Spin^c$-structure if and only if there exists an element $\omega\in H^2(X;\Z)$ such that $w_2(X)\equiv \omega \text{ mod }2$. The $Spin^c$-structures on $X$, if any exist, are parametrized by $2H^2(X;\Z)\oplus H^1(X;\Z_2)$. If $X$ is a complex manifold, it carries a canonical $Spin^c$-structure. In this case, the principal $U(1)$-bundle corresponds to the determinant line bundle of the complex tangent bundle of $X$.

If $X$ is a $Spin^c$-manifold, one can construct a complex spinor bundle $S(X)$ over it. The Levi-Civita connection on $X$ together with a connection $\nabla^c_X$ on $P_{U(1)}$ determine a connection on $S(X)$, which turns $S(X)$ into a Dirac bundle. There is an associated $Spin^c$-Dirac operator $D^c_X:\Gamma(S(X))\rightarrow \Gamma(S(X))$, where $ \Gamma(S(X))$ is the set of sections on $S(X)$. Given a Hermitian complex vector bundle $E$ over $X$, one also gets a twisted $Spin^c$-Dirac operator $D^c_{X,E}:\Gamma(S(X)\otimes E)\rightarrow \Gamma(S(X)\otimes E)$.

\subsection{The Atiyah-Patodi-Singer index theorem}

Eta-invariants were introduced in the context of the Atiyah-Patodi-Singer (APS) index theorems, where they appear as metric dependent boundary correction terms (see \cite{APSI75} and \cite{APSII75}).

Let $(W^{2n},g_W)$ be a compact $Spin^c$-manifold with boundary $\partial W=M$ such that the metric $g_W$ is of product form near the boundary and $g_M$ is the induced metric on $M$. Fix a connection $\nabla^c_W$ on the principal $U(1)$-bundle $P_{U(1)}\rightarrow W$ corresponding to the $Spin^c$-structure, which is constant in normal direction near the boundary and with curvature form $\Omega^c$. As $W$ is of even dimension, there is a direct sum splitting $S(W)=S^+(W)\oplus S^-(W)$. If we subject the restricted $Spin^c$-Dirac operator to APS-boundary conditions, we obtain an operator $D^{c,+}_W:\Gamma(S^+(W))\rightarrow \Gamma(S^-(W))$ and an adjoint operator $(D^{c,+}_W)^*$, both of which have finite dimensional kernel (see \cite[Appendix D]{LM89}).

Then the index $\text{ind} (D^{c,+}_W):=\text{dim(ker(}D^{c,+}_W))-\text{dim(ker(}(D^{c,+}_W)^*))\in \Z$ is well-defined, and by the APS index theorem (\cite[Theorem 4.2]{APSI75}), it is given by $$\text{ind} (D^{c,+}_W)=\int_We^{\frac{1}{2}c}\hat{\mathcal{A}}(W,g_W)-\frac{1}{2}\Big(h(M,g_M)+\eta(M,g_M)\Big).$$ Here $ \hat{\mathcal{A}}(W,g_W)$ is the $\hat{\mathcal{A}}$-series in Pontryagin forms, $c:=\frac{1}{2\pi}\Omega^c$ denotes the first Chern form in terms of the curvature form $\Omega^c$ of $\nabla^c$, $h(M,g_M):=\text{dim}(\text{ker}(D^c_M))$ and $\eta(M,g_M)$ is the eta-invariant of $D^c_M$. To introduce it, we need to consider $$\eta(z):=\sum_\lambda \frac{\text{sign}(\lambda)}{|\lambda|^z}, \qquad z\in\C, \,\text{Re}(z)\gg 0,$$ where the sum is running over all non-zero eigenvalues $\lambda$ of $D^c_M$. This function extends in a unique way to a meromorphic function (also denoted by $\eta$) which is holomorphic at $z=0$ (see \cite[\textsection 2]{APSI75}) and we may define $\eta(M,g_M):=\eta(0)$.

The twisted case is completely analogous. Let $E$ be a Hermitian complex vector bundle over $W$ and $\nabla^E$ a Hermitian connection on it, which is constant in normal direction near the boundary. The twisted APS index formula (\cite[(4.3)]{APSI75}) is then given by $$\text{ind}(D^{c,+}_{W,E})=\int_Wch(E,\nabla^E)e^{\frac{1}{2}c}\hat{\mathcal{A}}(W,g_W)-\frac{1}{2}\Big(h_E(M,g_M)+\eta_E(M,g_M)\Big).$$ Here $h_E(M,g_M):=\text{dim}(\text{ker}(D^c_{M,E}))$ is the kernel, $ch(E,\nabla^E):=\text{tr}(\text{exp}(\frac{i\Omega^c_E}{2\pi}))$ is the Chern character in terms of the curvature form $\Omega^c_E$ of $\nabla^E$ and $\eta_E(M,g_M)$ is defined as above with $D^{c}_{M,E}$.

We have the following well known vanishing theorem (see \cite{LM89} and \cite{APSII75}).

\begin{theorem}\label{thm:vanish}
Let $(W^{2n},g_W)$, $(M,g_M)$ and $E$ be as above, with $\Omega^c_W$ and $\Omega^c_M$ the curvature forms corresponding to the principal $U(1)$-bundle of the $Spin^c$-structure respectively. Suppose $E$ is flat. Then if $scal_{g_M}>0$ and $\Omega^c_M=0$, we have $h(M,g_M)=0$ and $h_E(M,g_M)=0$. If $scal_{g_W}>0$ and $\Omega^c_W=0$, then $\emph{ind}(D^{c,+}_{W})=0$ and $\emph{ind}(D^{c,+}_{W,E})=0$. 
\end{theorem}

Now suppose $M$ is connected. Let $\alpha:\pi_1(M)\rightarrow U(k)$ be a unitary representation, so that it determines a flat vector bundle $E_\alpha:=\tilde{M}\times_{\pi_1(M)}\C^k$ over $M$, where $\tilde{M}$ is the universal cover of $M$ (see \cite{K87} (1.2.4)). We can then define the reduced eta-invariant $$\tilde{\eta}_{E_\alpha}(M,g_M):=\eta_{E_\alpha}(M,g_M)-k\cdot\eta(M,g_M).$$

This invariant will allow us to tell path components apart in the space of Riemannian metrics of positive scalar curvature. Indeed, if $\mathcal{R}_{scal>0}(M)$ denotes the set of Riemannian metrics on $M$ with positive scalar curvature, we have the following result (see \cite[p. 417]{APSII75}, \cite[Proposition 3.3]{DGA19}).

\begin{proposition}\label{prop:etainvsame}
Let $M^{2n-1}$ be a closed connected $Spin^c$-manifold, $\nabla^c_M$ a flat connection on $P_{U(1)}\rightarrow M$ and $\alpha:\pi_1(M)\rightarrow U(k)$ a unitary representation. Let $g_0,g_1\in\mathcal{R}_{scal>0}(M)$ be in the same path component. Then $\tilde{\eta}_{\alpha}(M,g_0)=\tilde{\eta}_{\alpha}(M,g_1)$.
\end{proposition}

The reduced eta-invariant is in general difficult to compute. In an equivariant setting however, there are explicit formulas which allow us to determine them. Let $(M,g_M)$ be a closed Riemannian manifold and suppose a compact Lie group $G$ acts by isometries on it. Then one can define $G$-equivariant $Spin^c$-structures by requiring all of the corresponding bundles to be $G$-equivariant. If furthermore the unitary connection on the associated principal $U(1)$-bundle is $G$-equivariant, then one gets a $G$-equivariant spinor bundle $S_G(M)$ and a corresponding $G$-equivariant $Spin^c$-Dirac operator $D^{c,G}_M$. Fix $g\in G$ and let $\lambda$ be an eigenvalue of $D^{c,G}_M$. The action of $g$ on the Eigenspace $E_\lambda$ will be denoted by $g^\#_\lambda$. One can then consider $$\eta_g(z):=\sum_{0\neq\lambda \in E_\lambda}\frac{\text{sign}(\lambda)\cdot \text{tr}(g^\#_\lambda)}{|\lambda|^z}$$ which converges absolutely for $\text{Re}(z)\gg 0$ and define $\eta(M,g_M)_g:=\eta_g(0)$, the equivariant eta-invariant at $g$ (after meromorphic extension, see \cite{D78}).

Suppose now that $M$ is connected and that $G=\pi_1(M)$ is finite. The metric $g_M$ lifts to a metric $\tilde{g}_{\tilde{M}}$ on the universal cover $\tilde{M}$, as does the $Spin^c$-structre on $M$. Recall that $\alpha$ is a unitary representation and $E_\alpha$ the corresponding flat vector bundle, defined above. In this situation, we get the following formula (see \cite[Theorem 3.4.]{D78}). $$\eta_{E_\alpha}(M,g_M)=\frac{1}{|G|}\sum_{g\in G}\eta(\tilde{M},\tilde{g}_{\tilde{M}})_g\cdot \chi_\alpha(g),$$ where $\chi_\alpha:G\rightarrow \C:g\mapsto tr(\alpha(g))$ is the character of $\alpha$ and $|G|$ the order of $G$.

\subsection{Donnelly's Fixed Point Formula}

Let $(W^{2n},g_W)$ and $(M^{2n-1},g_M)$ be as in the previous paragraph and suppose a compact Lie group $G$ acts on $W$ by isometries and preserves the $Spin^c$-structure. Assume also that $g_W$ is of product form near the boundary $\partial W=M$. Let $D^{c,+}_W$ be the $G$-equivariant $Spin^c$-Dirac operator on $(W,g_W)$ subject to the APS-boundary condition and $D_M$ the $G$-equivariant $Spin^c$-Dirac operator on $(M,g_M)$ (with corresponding adjoint operators $(D^{c,+}_W)^*$ and $D^*_M$). In this case, the kernels of these operators are $G$-representations and the equivariant index $\text{ind}_G(D^{c,+}_W):=\text{ker}(D^{c,+}_W)-\text{ker}((D^{c,+}_W)^*)$ is a virtual complex $G$-representation in the ring of representations $R(G)$. For $g\in G$, let $h_g$ be the character of the $G$-representation $\text{ker}(D_M)$ at $g$ and $\eta(M,g_M)_g$ the equivariant eta-invariant of $D_M$ at $g$. Then we have the following fixed point formula (see \cite[Theorem 1.2.]{D78}).
$$ \text{ind}_G(D^{c,+}_W)_g + \frac{1}{2}\Big(h_g+\eta(M,g_M)_g\Big)=\sum_{N\subset W^g}a_N,$$ where $N$ denotes a connected component of the fixed point set $W^g$ and $a_N$ is a so-called local contribution at $N$. In our specific case, this local contribution can be computed explicitely (see \cite[Proposition 3.5]{DGA19}).

\begin{proposition}\label{prop:localcontr}
Let $(W^{2n},g_W)$ and $(M^{2n-1},g_M)$ be as above, and suppose $\tau\in G$ induces an involution on $W$ with isolated fixed points which lie in the interior of $W$. Then the local contribution at any fixed point $p$ is given by $a_{\{p\}}=2^{-n}$.
\end{proposition}

We can now compute the reduced eta-invariant in the case of the Brieskorn quotients. Let $d\in\N$ be even and $\epsilon\neq 0$.

Let $W^6_{\epsilon}(d)$ be equipped with a $G$-invariant metric $g_W$, where $G$ is a closed subgroup of $\Z_{2d}\times O(3)$ containing the involution $\tau$ defined above. Fix a Hermitian metric on its tangent bundle which is $G$-invariant. Consider the $G$-equivariant $Spin^c$-structure induced by its complex manifold structure and the Hermitian metric. The induced $G$-equivariant $Spin^c$-structure on the boundary $\Sigma^5_{\epsilon}(d)=\partial W^6_{\epsilon}(d)$ descends to a $Spin^c$-structure on $N^5_{\epsilon}(d)=\Sigma^5_{\epsilon}(d)/\tau$, which we will call its preferred $Spin^c$-structure.

Now fix a flat connection on the principal $U(1)$-bundles over $W^6_{\epsilon}(d)$, $\Sigma^5_{\epsilon}(d)$ and $N^5_{\epsilon}(d)$ corresponding to the $Spin^c$-structure (which we can do since the second cohomology groups with real coefficients are trivial and therefore the first Chern classes of the bundles vanish). If $\tau$ acts as $-Id$ on $\C$, then $E_\alpha:=\Sigma^5_{\epsilon}(d)\times_{<\tau>}\C$ is a non-trivial flat complex line bundle.

\begin{proposition}\label{prop:etainveps}
Let $0<\epsilon_0 < 1$ and $g_N$ be a metric of $scal>0$ on $N^5_\epsilon(d)$. Suppose that the lift $g_\Sigma$ of $g_N$ to $\Sigma^5_\epsilon(d)$ extends to a metric $g_W$ on $W^6_\epsilon(d)$ which is $\tau$-invariant, of product form near the boundary and of $scal>0$. Then $$\tilde{\eta}_{E_\alpha}(N^5_\epsilon(d),g_N)=-\frac{d}{4}.$$

\begin{proof}
Using the formulas from the previous section, we get for $G=\Z_2=\{1,\tau\}$:

$$\eta_{E_\alpha}(N^5_\epsilon(d),g_N)=\frac{1}{2}\Big(\eta(\Sigma^5_\epsilon(d),g_\Sigma)_1\cdot \chi_\alpha(1)+\eta(\Sigma^5_\epsilon(d),g_\Sigma)_\tau\cdot \chi_\alpha(\tau)\Big)$$
$$=\frac{1}{2}\Big(\eta(\Sigma^5_\epsilon(d),g_\Sigma)-\eta(\Sigma^5_\epsilon(d),g_\Sigma)_\tau\Big)$$
and using the trivial representation (whose character is always 1)
$$\eta(N^5_\epsilon(d),g_N)=\frac{1}{2}\Big(\eta(\Sigma^5_\epsilon(d),g_\Sigma)+\eta(\Sigma^5_\epsilon(d),g_\Sigma)_\tau\Big).$$
Hence we get $$\tilde{\eta}_{E_\alpha}(N^5_\epsilon(d),g_N)=\eta_{E_\alpha}(N^5_\epsilon(d),g_N)-1\cdot \eta(N^5_\epsilon(d),g_N)=-\eta(\Sigma^5_\epsilon(d),g_\Sigma)_\tau$$
Now since the metrics on $\Sigma^5_\epsilon(d)$ and $W^6_\epsilon(d)$ both are $\tau$-invariant and of $scal>0$, applying Theorem \ref{thm:vanish}, we see that the equivariant index and kernel vanish. Using the fixed point formula from above we get $\eta(\Sigma^5_\epsilon(d))_\tau=2\sum_{i=1}^d a_{\{p_i\}}$, where the $p_i$ are the isolated fixed points of the action of $\tau$ on $W^6_\epsilon(d)$. Finally, using Proposition \ref{prop:localcontr} to compute $a_{\{p_i\}}=2^{-n}$ with $n=3$, we get the desired result.

\end{proof}

\end{proposition}

\section{Moduli space of Riemannian metrics}

For a general treatment on moduli spaces, see \cite{TW15}.

Let $M$ be a closed smooth manifold and $\mathcal{R}(M)$ the set of Riemannian metrics on $M$, equipped with the $C^\infty$-topology of uniform convergence of all the derivatives. By taking pull-backs of metrics, we get an action of the group of diffeomorphisms Diff($M$) on $\mathcal{R}(M)$. The quotient space $\mathcal{M}(M):=\mathcal{R}(M)/\text{Diff}(M)$ is called the \textit{moduli space} of $M$. If we restrict to metrics of nonnegative sectional curvature or positive Ricci curvature, we obtain $\mathcal{R}_{\text{sec}\geq 0}(M)$ and $\mathcal{R}_{\text{Ric}> 0}(M)$ respectively. Restricting on isometry classes of metrics of nonnegative sectional curvature or positive Ricci curvature, we obtain corresponding moduli spaces $\mathcal{M}_{\text{sec}\geq 0}(M)$ and $\mathcal{M}_{\text{Ric}> 0}(M)$.

For technical reasons, we will actually need to work in the quotient space $\mathcal{R}_{sec\geq 0}(M)/\mathcal{D}$, where $\mathcal{D}\subset \text{Diff}(M)$ is the subgroup of diffeomorphisms which preserve the $Spin^c$-structure on $M$.

First we need to determine the reduced eta-invariant of $(N^5_0(d),g^d)$, where $g^d$ is the modified Grove-Ziller metric from \textsection \ref{sec:metrics}. Since $H^2(N^5_0(d);\Z)\cong \Z_2$, there exists a unique flat non-trivial complex line bundle over $N^5_0(d)$ which we will denote by $\alpha$. Fix $0<\epsilon_0<1$.

\begin{lemma}\label{lemma:etainvcomp}
Let $(N^5_0(d),g^d)$ be equipped with the preferred $Spin^c$-structure and a flat connection on the associated principal $U(1)$-bundle. Then $$\tilde{\eta}_\alpha(N^5_0(d),g^d)=-\frac{d}{4}.$$

\begin{proof}
Let $$\tilde{Z}(d)=\bigcup_{0\leq \epsilon \leq \epsilon_0}N^5_\epsilon(d)$$ be equipped with the extended metric $h$ as in \textsection \ref{sec:metrics}. Let $g^d_t$ and $h_t$ denote the Cheeger deformations of $g^d$ and $h$ respectively and $t_0$ the parameter from Lemma \ref{lemma:metrics}. We construct a path in $\mathcal{R}_{scal>0}(N^5_{\epsilon_0}(d))$, connecting $(N^5_0(d),g^d)$ to $(N^5_{\epsilon_0}(d),h_{t_0}|_{N^5_{\epsilon_0}(d)})$, for which we can actually compute the reduced eta-invariant as in Proposition \ref{prop:etainveps}. Let $\phi_\epsilon:N^5_{\epsilon_0}(d)\rightarrow N^5_{\epsilon}(d)$ be a smooth family of diffeomorphisms for $0<\epsilon<\epsilon_0$.

Using Lemma \ref{lemma:metrics}, we can now define two paths of metrics of positive scalar curvature on $N^5_{\epsilon_0}(d)$ in the following way.

\begin{align*}
\tilde{\gamma}_1: & \text{ }[0,t_0]\rightarrow \mathcal{R}_{scal>0}(N^5_{\epsilon_0}(d)) \\ & \quad \text{ } t \mapsto \phi^*_0(g^d_t)
\end{align*}
and
\begin{align*}
\tilde{\gamma}_2: & \text{ } [0,\epsilon_0]\rightarrow \mathcal{R}_{scal>0}(N^5_{\epsilon_0}(d)) \\ & \quad \text{ } \epsilon \mapsto \phi^*_\epsilon(h_{t_0}|_{N^5_{\epsilon}(d)})
\end{align*}

\noindent Since $\tilde{\gamma}_1(t_0)=\phi^*_0(g^d_{t_0})=\phi^*_0(h_{t_0}|_{N^5_{0}(d)}))=\tilde{\gamma}_2(0)$, we can concatenate these two paths to get a path $\tilde{\gamma}$ in $\mathcal{R}_{scal>0}(N^5_{\epsilon_0}(d))$ with endpoints $\phi^*_0(g^d)$ and $\tilde{h}_{t_0}|_{N^5_{\epsilon_0}(d)}$. By Lemma \ref{lemma:metrics} and the discussion above it, $h_{t_0}|_{N^5_{\epsilon_0}(d)}$ lifts to a metric on $\Sigma^5_{\epsilon_0}(d)$ which extends to a metric on $W^6_{\epsilon_0}$ satisfying the conditions of Proposition \ref{prop:etainveps}. By Propostion \ref{prop:etainvsame} and using the fact that the reduced eta-invariant is preserved under pull-backs by diffeomorphisms, we finally get:

$$\tilde{\eta}_\alpha(N^5_0(d),g^d)=\tilde{\eta}_\alpha(N^5_{\epsilon_0}(d),\phi^*_0(g^d))=\tilde{\eta}_\alpha(N^5_{\epsilon_0}(d),h_{t_0}|_{N^5_{\epsilon_0}(d)})=-\frac{d}{4}.$$
\end{proof}

\end{lemma}

We can now prove our main theorem.

Fix $d\in\{0,2,4,6,8\}$, $d_0\neq d_1\in \N$ such that $d_i\equiv \pm d \text{ mod } 16$. Therefore, by Theorem \ref{thm:classification5mfds} and Lemma \ref{lemma:diffquo}, there exist diffeomorphisms $\psi_i:N^5_0(d)\rightarrow N^5_0(d_i)$ for $i=0,1$. Denote by $h_i=\psi^*_i(g^{d_i})$ the pull-back metric for $i=0,1$. Let $\mathcal{D}\subset \text{Diff}(N^5_0(d))$ denote the subgroup of diffeomorphisms preserving the $Spin^c$-structure on $N^5_0(d)$.

Since $\mathcal{D}$ has finite index in Diff$(N^5_0(d))$, if $\mathcal{R}_{sec\geq 0}(N^5_0(d))/\mathcal{D}$ has infinitely many path components, it will follow that $\mathcal{M}_{\text{sec}\geq 0}(N^5_0(d))$ has infinitely many path components as well.

Assume now that there is a path $\gamma:[0,1]\rightarrow \mathcal{R}_{sec\geq 0}(N^5_0(d))/\mathcal{D}$ connecting $[h_0]$ and $[h_1]$.  As a consequence of Ebin's slice theorem, this path can be lifted to a path $\tilde{\gamma}$ in $\mathcal{R}_{sec\geq 0}(N^5_0(d))$ such that $\tilde{\gamma}(0)=h_0$ and $\tilde{\gamma}(0)=\phi^*(h_1)$ for some $\phi\in\mathcal{D}$ (see \cite{E70} and \cite[Proposition 4.6]{CK19}). By B\"{o}hm and Wilking, the path $\tilde{\gamma}$ evolves instantly to a path in $\mathcal{R}_{Ric> 0}(N^5_0(d))$ under the Ricci flow (see \cite[Theorem A]{BW07}). If we concatenate this resulting path with the orbits of the endpoints of $\tilde{\gamma}$, we obtain a path $\gamma'$ in $\mathcal{R}_{scal> 0}(N^5_0(d))$ with endpoints $\gamma'(0)=h_0$ and $\gamma'(1)=\phi^*(h_1)$ as well.

Using Lemma \ref{lemma:etainvcomp}, we can now compute the reduced eta-invariants of the endpoints of $\gamma'$. We immediately see that $\tilde{\eta}_\alpha(N^5_0(d_0),h_0)=-\frac{d_0}{4}\neq-\frac{d_1}{4}=\tilde{\eta}_\alpha(N^5_0(d_1),\phi^*(h_1))$. Hence $h_0$ and $\phi^*(h_1)$ lie in different path components in $\mathcal{R}_{scal>0}(N^5_{0}(d))$, which contradicts the construction we have obtained from our initial assumption.

By Theorem \ref{thm:classification5mfds}, Lemma \ref{lemma:diffquo} and the above remark on the link between path components of $\mathcal{R}_{sec\geq 0}(N^5_0(d))/\mathcal{D}$ and $\mathcal{M}_{\text{sec}\geq 0}(N^5_0(d))$, it follows that the moduli space $\mathcal{M}_{\text{sec}\geq 0}(N^5_0(d))$ has infinitely many path components. If we focus on the evolved path under the Ricci flow, we see that the same argument applies to $\mathcal{R}_{Ric > 0}(N^5_0(d))/\mathcal{D}$ and hence $\mathcal{M}_{Ric > 0}(N^5_0(d))$ has infinitely many path components as well. \proofend

\noindent\textsc{Department of Mathematics, University of Fribourg, Switzerland}
\emph{E-mail address: } jonathan.wermelinger@unifr.ch


\begin{thebibliography}{9}

\bibitem[AB15]{AB15} Alexandrino, Marcos M.; Bettiol, Renato G. Lie groups and geometric aspects of isometric actions. Springer, Cham, 2015. x+213 pp. ISBN: 978-3-319-16612-4; 978-3-319-16613-1 MR3362465

\bibitem[ASIII68]{ASIII68} Atiyah, M. F.; Singer, I. M. The Index of Elliptic Operators: III. Annals of Mathematics, Second Series, Vol. 87, No. 3 (May, 1968), pp. 546-604 .

\bibitem[APSI75]{APSI75} Atiyah, M. F.; Patodi, V. K.; Singer, I. M. Spectral asymmetry and Riemannian geometry. I. Math. Proc. Cambridge Philos. Soc. 77 (1975), 43–69

\bibitem[APSII75]{APSII75} Atiyah, M. F.; Patodi, V. K.; Singer, I. M. Spectral asymmetry and Riemannian geometry. II. Math. Proc. Cambridge Philos. Soc. 78 (1975), no. 3, 405–432.


\bibitem[BW07]{BW07} B\"ohm, C. and Wilking, B. Nonnegatively curved manifolds with finite fundamental groups admit metrics with positive Ricci curvature, Geom. Funct. Anal. 17 (2007) 665–681

\bibitem[CK19]{CK19} Corro, Diego; Kordass, Jan-Bernhard. Short survey on the existence of slices for the space of Riemannian metrics. Preprint, arxiv:1904.07031,(2019).

\bibitem[DGA19]{DGA19} Dessai, Anand; Gonz\'{a}lez-\'{A}lvaro, David. Moduli space of metrics of nonnegative sectional or positive Ricci curvature on homotopy real projective spaces. Preprint, arXiv:1902.08919 (2019)

\bibitem[D78]{D78} Donnelly, Harold. Eta invariants for $G$-spaces. Indiana Univ. Math. J. 27 (1978), no. 6, 889--918. MR0511246

\bibitem[E70]{E70} Ebin, David G. The manifold of Riemannian metrics. 1970 Global Analysis (Proc. Sympos. Pure Math., Vol. XV, Berkeley, Calif., 1968) pp. 11--40 Amer. Math. Soc., Providence, R.I. MR0267604 

\bibitem[GT98]{GT98} Geiges, H.; Thomas, C.B. Contact Topology and the Structure of 5-Manifolds with $\pi_1=\Z_2$, Ann. Inst. Fourier, Grenoble, 48, 4 (1998), 1167-1188

\bibitem[G20]{G20} Goodman, McFeely Jackson. Moduli spaces of Ricci positive metrics in dimension five. Preprint, arXiv:2002.00333 (2020)

\bibitem[GVWZ06]{GVWZ06} Grove, Karsten; Verdiani, Luigi; Wilking, Burkhard; Ziller, Wolfgang. Non-negative curvature obstructions in cohomogeneity one and the Kervaire spheres. Ann. Sc. Norm. Super. Pisa Cl. Sci. (5) 5 (2006), no. 2, 159-170. MR2244696

\bibitem[GZ00]{GZ00} Grove, Karsten; Ziller, Wolfgang. Curvature and symmetry of Milnor spheres. Ann. of Math. (2) 152 (2000), no. 1, 331--367. MR1792298



\bibitem[KT90]{KT90} R. C. Kirby; L. R. Taylor. Pin structures on low-dimensional manifolds. Geometry of low-dimensional manifolds, 2 (Durham, 1989), London Math. Soc. Lecture Note Ser., vol. 151, Cambridge Univ. Press, Cambridge, 1990, pp. 177–242

\bibitem[K87]{K87} Kobayashi, Shoshichi. Differential Geometry of Complex Vector Bundles. Princeton University Press, Princeton Legacy Library 783 (1987)

\bibitem[HM68]{HM68} Hirzebruch, F.; Mayer, K.H. $O(n)$-Mannigfaltigkeiten, exotische Sphären und Singularitäten. Lecture Notes in Mathematics 57, Berlin-Heidelberg-New York: Springer 1968.

\bibitem[LM89]{LM89} Lawson, H. Blaine, Jr.; Michelsohn, Marie-Louise. Spin geometry. Princeton Mathematical Series, 38. Princeton University Press, Princeton, NJ, 1989. {\rm xii}+427 pp. ISBN: 0-691-08542-0 MR1031992 

\bibitem[M68]{M68} Milnor, John. Singular points of complex hypersurfaces. Annals of Mathematics Studies, No. 61 Princeton University Press, Princeton, N.J.; University of Tokyo Press, Tokyo 1968.

\bibitem[S61]{S61} Smale, Stephen. Generalized Poincar\'{e}'s conjecture in dimensions greater than four. Ann. of Math. (2). 74 (2): 391406 (1961).

\bibitem[S12]{S12} Su, Yang. Free involutions on $S^2\times S^3$. Geom Dedicata (2012) 159:11-28.

\bibitem[TWi17]{TWi17} Tuschmann, Wilderich; Wiemeler, Michael. On the topology of moduli spaces of non-negatively curved Riemannian metrics. Preprint, arXiv:1712.07052 (2017)

\bibitem[TW15]{TW15} Tuschmann, Wilderich; Wraith, David J. Moduli spaces of Riemannian metrics. Second corrected printing. Oberwolfach Seminars, 46. Birkh\"auser Verlag, Basel, 2015.



\end{thebibliography}
\end{document}